\documentclass[12pt]{article}
\author{Dave Platt\footnote{Supported by Australian RC Discovery Project DP160100932 and  EPSRC Grant EP/K034383/1.}\\ School of Mathematics, \\ University of Bristol, Bristol, UK\\dave.platt@bris.ac.uk\and  Tim Trudgian$^{*}$\footnote{Supported by Australian RC Future Fellowship FT160100094.} \\
School of Science\\ The University of New South Wales Canberra, Australia \\
  t.trudgian@adfa.edu.au}
  \title{The Riemann hypothesis is true up to $3\cdot 10^{12}$}
\usepackage{indentfirst}
\usepackage{url}
\usepackage{enumerate}
\usepackage{amsthm}
\usepackage{amsmath}
\usepackage{comment}
\usepackage{fullpage}
\usepackage{amssymb}
\usepackage{booktabs}

\newtheorem{thm}{Theorem}

\newtheorem{cor}{Corollary}

\DeclareMathOperator{\li}{li}
\newlength{\bibitemsep}
\newlength{\bibparskip}
\let\oldthebibliography\thebibliography
\renewcommand\thebibliography[1]{%
  \oldthebibliography{#1}%
  \setlength{\parskip}{\bibitemsep}%
  \setlength{\itemsep}{\bibparskip}%
}

\begin{document}
\maketitle
\begin{abstract}
\noindent 
We verify numerically, in a rigorous way using interval arithmetic, that the Riemann hypothesis is true up to height $3\cdot10^{12}$. That is, all zeroes $\beta + i\gamma$ of the Riemann zeta-function with $0<\gamma\leq 3\cdot 10^{12}$ have $\beta = 1/2$.
\end{abstract}
\section{Introduction}
\noindent
The Riemann zeta-function $\zeta(s)$ has trivial zeroes at $s=-2,-4,-6\ldots$, and non-trivial zeroes in the strip $0< \sigma <1$, where here, and hereafter $s= \sigma + it$. The Riemann hypothesis asserts that all non-trivial zeroes  $\rho = \beta + i\gamma$ have $\beta = 1/2$. In the absence of a proof, it is extremely important to obtain partial verifications of the Riemann hypothesis. To that end, define $H$ as the largest number for which it is known that all zeroes $\rho = \beta + i\gamma$ with $0<\gamma \leq H$ have $\beta = 1/2$.
This problem has a long history; for a glimpse of this we refer the reader to \cite[p.\ 2]{Gourdon}.

In recent years, three calculations have been referenced frequently in the literature. The first result is by Wedeniwski\footnote{It is difficult to pinpoint the height claimed in these calculations. The third slide of \cite{Wedeniwski}, dated $20$th May $2003$, claims the first $200$ billion zeroes were checked, equating to $H=5.72\ldots\cdot 10^{10}$. The end of \cite{Wedeniwski} contains the statement that the project computed $385$ billion zeroes or $H=1.07\ldots\cdot 10^{11}$ if we assume these were the lowest lying zeroes. A later version of the same slides from $12$th November $2003$ mentions $561$ billion zeroes or $H=1.53\ldots\cdot 10^{11}$. To add to the confusion, the link \cite{Wiki}  gives $900$ billion zeroes and this is leads to the quoted $H=2.41\ldots\cdot 10^{11}$.} \cite{Wedeniwski} in 2004, with $H=2.41\ldots\cdot 10^{11}$. 
The second is by Gourdon \cite{Gourdon}, also in 2004, which establishes $H= 2.44\ldots \cdot 10^{12}$. The third is by the first author in 2017 \cite{Platt} (see also \cite{PlattThesis}), with $H= 3.06\ldots\cdot 10^{10}$.

Whilst the latest and lowest value for $H$ may appear to be a retrograde step, the computations in \cite{Platt} utilised interval arithmetic and rigorously derived truncation bounds to ensure the results claimed are correct. The two earlier results have the disadvantage that neither has been published in a peer reviewed journal. Furthermore, it is not clear how the computations were set up to avoid problematic accumulation of rounding and truncation errors. This concern was noted in works by Tao \cite{Tao} (see remarks after Theorem 1.5) and Helfgott \cite{Helfgott}. We also mention the result claimed by Franke et al. \cite{Frank} with $H= 10^{11}$: again, few computational details are given.

We now state our main theorem, which surpasses all aforementioned results.
\begin{thm}\label{bank}
The Riemann hypothesis is true up to height $3\,000\,175\,332\,800$. That is, the lowest $12\,363\,153\,437\,138$ non-trivial zeroes $\rho$ have $\Re\rho = 1/2$.
\end{thm}

We note that this independently verifies the results of Gourdon, Wedeniwski and Franke et al.\ and indeed goes $22\%$ higher than the largest of these. 

We have endeavoured to make this paper as short as possible. In \S \ref{plank} we outline some of the computational aspects underpinning Theorem \ref{bank}. In \S \ref{dank} we mention some results that are improved instantly with Theorem \ref{bank}.

\section{Theory and computation}\label{plank}
We used the algorithm described in \cite{Platt}. In common with all modern partial verifications of the Riemann hypothesis, the algorithm computes values of the completed zeta function on the half line and counts sign changes therein. Each sign change represents a zero of zeta on the half line. Using a variation of Turing's method (see \cite{Turing1,Turing2}), we can confirm that all the expected zeroes have been accounted for, so none lie off the half line and the Riemann hypothesis holds in the given range.

We rewrote the original code to utilise Arb \cite{ARB} in place of MPFI \cite{Revol2002} for two reasons. First, Arb is being actively maintained whereas MPFI is not. Second, Arb uses ball arithmetic in place of full interval arithmetic whence there is a space saving of roughly $50\%$, which make applications more cache friendly.

The other main change to the code that was used to reach $3\cdot 10^{10}$ in \cite{Platt} was that we made no attempt to isolate zeroes to any more precision than was absolutely necessary. A key motivation of \cite{Platt} was to generate a database of rigorously isolated zeroes to high precision, but to do so here would have added to the run time and, in any case, we had nowhere to store that many zeroes. Rather, once we had found a sign change in the completed zeta function indicating the presence of a zero of zeta on the half line, we did not use the machinery of Shannon--Whittaker--Nyquist to ``zoom in'' on that zero, we merely counted it and moved on. In fact, the default lattice sampling rate we used (about $0.01$) was sufficient to isolate $999\,997.5$ out of every $1\,000\,000$ zeroes.

All computations were performed on the University of Bristol's Bluecrystal Phase III \cite{ACRC} and the National Computing Infrastructure's Raijin and Gadi \cite{NCI} clusters. Both Raijin and Gadi nodes support Hyper-Threading meaning there are two logical processors per core. Since both of these logical devices share the same physical execution resources, one would not expect to see a two times speed up: we found we benefitted to the tune of about $15\%$. Bluecrystal does not have Hyper-Threading enabled. We also incorporated a minor improvement to our bound for $|\Gamma((\sigma+it)/2)|\exp(\pi t/4)$ (A.2 in \cite{Platt}) and spent some time optimising the computational parameters to work better $100$ times higher up the half line.  

In total we used some $7.5$ million core hours on $3.6$GHz Intel\textsuperscript{\textregistered} Xeon\textsuperscript{\textregistered} processors, so each GHz-hour processed a piece of the half line of length about $110\,000$. For comparison, Wedeniwski \cite{Wedeniwski} reports that the isolation of $561$ billion zeros took the equivalent of $2\,304$ years on $2$GHz Pentium\textsuperscript{\textregistered} $4$, so about $3\,800$ of the half line per GHz-hour whereas Gourdon's computation took the equivalent of $525$ days on a single $2.4$GHz Pentium\textsuperscript{\textregistered} $4$, so about $80\,000\,000$ of the half line per GHz-hour. The difference between the Wedeniwski computation and Gourdon's and our's shows the power of FFT based algorithms up against vanilla Riemann-Siegel. The fact that Gourdon's computation was $725$ times quicker than ours is down to the higher sampling rate we used ($25$ per zero in place of $1.2$) and the cost of multi precision rigorous numerics compared to hardware floating point. 

\section{Some instant wins}\label{dank}
Where researchers have used Gourdon's or Platt's $H$ as their starting point, this independent verification adds weight to their results and gives some explicit improvements.
In this section we include some results that are improved either instantly, or at least fairly easily, in light of Theorem \ref{bank}. We have not endeavoured to furnish an exhaustive list.

\subsection{Bounds on primes}

It is useful to have explicit estimates on the error term in the prime number theorem. Define $\psi(x) = \sum_{p^{m}\leq x} \Lambda(n)$, where $\Lambda(n)$ is the von Mangoldt function, and $\theta(x) = \sum_{p\leq x} \log p$. Rosser and Schoenfeld \cite{Ross} instigated a program of research to bound $\psi(x) -x$ explicitly. The current best results are due to Dusart \cite{Dusart} and Faber and Kadiri \cite{Faber_fix} for small values of $x$, Broadbent et al.\  \cite{Broadbent} for intermediate values,  and the authors \cite{PTnow} for larger values. These results were based on Gourdon's value of $H$ which is confirmed by and can be improved improved slightly with our Theorem \ref{bank}. 

Related to these bounds are Bertrand-type estimates: exhibiting a prime in intervals of the form $(x, x+ cx]$ for some $c$ and for all $x\geq x_{0}(c)$. The best results of this type are by Kadiri and Lumley \cite{KL}, and can now be improved with Theorem \ref{bank}.

Of course, on the Riemann hypothesis much more is known: Schoenfeld \cite{Sock} proved that
\begin{equation}\label{drank}
|\psi(x) -x| \leq \frac{1}{8\pi} x^{1/2} \log^{2} x, \quad (x > 59).
\end{equation}

B\"{u}the \cite[\S 7]{Booth} showed that if the Riemann hypothesis holds up to height $H$ then (\ref{drank}) holds for those $x$ such that $4.92\sqrt{x/\log x} \leq H$. Given this we may note a quick corollary.
\begin{cor}
The following bounds hold in the range indicated
\begin{equation*}
\begin{split}
|\psi(x) - x| &\leq \frac{\sqrt{x}}{8\pi} \log^{2} x, \quad (59<x\leq2.169\cdot 10^{25}),\\
|\theta - x| &\leq \frac{\sqrt{x}}{8\pi} \log^{2} x, \quad (599<x\leq2.169\cdot 10^{25}),\\
|\pi(x) - \li(x)| &\leq \frac{\sqrt{x}}{8\pi} \log^{2} x, \quad (2657<x\leq2.169\cdot 10^{25}).\\
\end{split}
\end{equation*}
\end{cor}

We also mention that some results from the expansive article by Bennett et al.\ \cite{Bennett} can be improved, as can bounds on weighted sums of $\Lambda(n)$ as given by Ramar\'{e} \cite{Ram}.

\subsection{Zero-free regions and zero density estimates}
%
It is known that there are no zeroes $\rho= \beta + i\gamma$ in the region $\beta \geq 1 - \frac{1}{R \log \gamma}$ for all $\gamma>3$. The current best value is $R = 5.573412$ by Mossinghoff and Trudgian \cite{Hoff}. This can be improved slightly with Theorem \ref{bank} and will be tackled in a future paper of Mossinghoff and the second author, along with the explicit version of the Vinogradov--Korobov zero-free region, which was proved by Ford \cite{Ford}.

Not only is the zeta-function non-zero in regions close to $\sigma =1$, its reciprocal does not grow too quickly. This has been quantified, explicitly, by the second author in \cite{Log}. These results depend not just on the value of $H$, and hence are improvable by Theorem \ref{bank}, but also on the size of the zero-free constant $R$ mentioned above.

Let $N(\sigma, T)$ count the number of zeroes with $\beta > \sigma$ and $0<\gamma\leq T$. Explicit estimates on $N(\sigma, T)$ have been given by Kadiri \cite{KS}, Kadiri, Lumley, and Ng \cite{KLN}, and Simoni\v{c} \cite{Aleks}. All of these results could be improved with Theorem \ref{bank}.

\subsection{Oscillations in the prime number theorem}
The prime number theorem gives $\pi(x) \sim \li(x)$ and $\theta(x) \sim x$. Littlewood \cite{Littlewood} showed that there are infinitely many sign changes in the differences $\pi(x) - \li(x)$ and $\theta(x) -x$. The history of estimating the first sign change of the first of these differences is rich: see \cite{Saouter,Smith} for further details. These results, as well as the results of the authors \cite{PT} for the first sign change of $\theta(x) -x$ could potentially be improved with Theorem \ref{bank}. Such an improvement would only be meagre, though, owing to the known small gap between unconditional results and those contingent on the Riemann hypothesis.

The aforementioned results show that the first sign changes cannot be too large. Numerical work shows that the first sign changes cannot be too small, either. B\"{u}the in \cite{Booth2} shows that $\theta(x) <x$ and $\pi(x) < \li(x)$ for $2\leq x\leq10^{19}$. We note that these results relied on the value of $H=10^{11}$ in Franke et al.\ \cite{Frank} and so, could also be improved with our Theorem \ref{bank}.
\subsection{The de Bruijn--Newman constant}
For $t\in \mathbb{R}$, let
\begin{equation*}
H_{t}(z)= \int_{0}^{\infty} e^{t u^{2}} \Phi(u) \cos(zu)\, du, \quad  \Phi(u) = \sum_{n=1}^{\infty} \left( 2 \pi^{2} n^{4} e^{9u} - 3\pi n^{2} e^{5u}\right) \exp\left(-\pi n^{2} e^{4u}\right).
\end{equation*}
The de Bruijn--Newman constant $\Lambda$ is the real number for which all the zeroes of $H_{t}$ are real when $t\geq \Lambda$. 
The Riemann hypothesis is equivalent to the statement that $\Lambda \leq 0$. Rodgers and Tao \cite{RT} proved that $\Lambda \geq 0$. A history of  bounds on $\Lambda$ is given in \cite{RT} and \cite{Poly}. 

The 15th Polymath Project \cite{Poly} contains some calculations with the de Bruijn--Newman constant: the authors prove that $\Lambda \leq 0.22$. We note that we can make an instant, but very mild, improvement on this. The second row in Table 1 on page 65 of \cite{Poly} shows\footnote{Note that the $X$ in Table 1 corresponds to $2H$. } that one may take $\Lambda \leq 0.2$ provided one has shown $H> 2.51\cdot 10^{12}$. This leads to the following.
\begin{cor}
We have $\Lambda \leq 0.2.$
\end{cor}

The next entry in Table 1 of \cite{Poly} is conditional on taking $H$ a little higher than $10^{13}$, which of course, is not achieved by Theorem \ref{bank}. This would enable one to prove $\Lambda < 0.19$. Given that our value of $H$ falls between the entries in this table, it is possible that some extra decimals could be wrought out of the calculation. We have not pursued this.

\subsubsection*{Acknowledgements}

The authors would like to thank the technical staff at NCI and the University of Bristol ACRC for their invaluable support and advice. This mammoth computation would not have been possible without the allocation of machine hours on Raijin and Gadi via NCMAS, INTERSECT, and the UNSW Resource Allocation Scheme.


 \end{document}